\documentclass[a4paper,10pt]{article}
\usepackage{amssymb, amsmath, amsthm, latexsym, verbatim, enumerate}

\newcommand \G{\Gamma}

\newcommand \la{\lambda}
\newcommand \id{\mathrm{id}}

\newcommand \Ker{\mathrm{Ker}}
\newcommand \Ric{\mathrm{Ric}}

\newcommand \Span{\mathrm{Span}}

\newcommand \Tr{\mathrm{Tr}}

\newcommand \diag{\mathrm{diag}}

\renewcommand \a{\alpha}
\renewcommand \b{\beta}

\newcommand \<{\langle}
\renewcommand \>{\rangle}
\newcommand \g{\mathfrak{g}}
\newcommand \n{\mathfrak{n}}
\newcommand \tn{\tilde{\mathfrak{n}}}
\newcommand \ad{\mathrm{ad}}
\newcommand \z{\mathfrak{z}}

\theoremstyle{plain}

\newtheorem*{theorem*}{Theorem}
\newtheorem*{corollary*}{Corollary}
\newtheorem{lemma}{Lemma}

\newtheorem*{proposition*}{Proposition}

\newtheorem*{namedtheorem}{\theoremname}
\newcommand{\theoremname}{te}

\theoremstyle{definition}

\newtheorem*{definition*}{Definition}

\theoremstyle{remark}

\begin{document}

\title{Harmonic homogeneous manifolds of nonpositive curvature
\footnote{2000 \emph{Mathematics Subject Classification.} Primary: 53C30, 53C25.}}

\author{Y.Nikolayevsky\thanks{Work supported by the ARC Discovery grant DP0342758.}}

\date{}

\maketitle

\begin{abstract} A Riemannian manifold is called \emph{harmonic} if its volume density
function expressed in polar coordinates centered at any point of the manifold is radial.
Flat and rank-one symmetric spaces are harmonic. The converse
(the \emph{Lichnerowicz Conjecture}) is true for manifolds of nonnegative scalar curvature
and for some other classes of manifolds, but is not true in general: there exists a
family of homogeneous harmonic spaces, the Damek-Ricci spaces, containing noncompact
rank-one symmetric spaces, as well as infinitely many nonsymmetric examples.
We prove that a harmonic homogeneous manifold of nonpositive curvature is either flat,
or is isometric to a Damek-Ricci space.
\end{abstract}

\section{Introduction}\label{s:intro}

Let $M^n$ be a Riemannian manifold, with $\nabla$ the Levi-Civita connection and $R$ the
curvature tensor. Let $\gamma(t)$ be a unit speed geodesic of $M^n$, with
$\gamma(0) = x \in M^n, \; \dot{\gamma}(0) = X \in T_xM^n$. Define a \emph{Lagrange tensor}
$A(t) \in \mathrm{End}(T_{\gamma(t)}M^n \cap \dot\gamma(t)^\perp)$ along $\gamma$ by
$$
\nabla_{\dot\gamma(t)}\nabla_{\dot\gamma(t)} A(t) + R_{\dot\gamma(t)} \circ A =0,\quad A(0) = 0, \quad \nabla_XA = \id_{X^\perp},
$$
where $R_{\dot\gamma(t)} = R(\cdot, \dot\gamma(t))\dot\gamma(t)$ is the Jacobi operator.
The function $V_{x, X}(t) = \det A(t)$ is the \emph{volume density function}.

\begin{definition*}
A Riemannian manifold $M^n$ is called \emph{harmonic} if its volume density functions
$V_{x, X}(t)$ depends neither of $x$, nor of $X$.
\end{definition*}

Equivalently, for any point $x \in M^n$, there exists a nonconstant harmonic function
defined on a punctured neighborhood of $x$ and depending only on the distance to $x$;
the mean (or the scalar) curvature of small geodesic spheres depends only on the radius
(see \cite[Ch.2.6]{BTV}, \cite[Ch.6]{B} for other equivalent definitions).

Expanding the volume density functions in the Taylor series one gets an infinite sequence
of conditions, the \emph{Ledger formulas}, on the curvature tensor and its covariant derivatives,
first two of which being
\begin{equation}\label{eq:ledger}
\Ric (X,X) = \Tr R_X = C \, \|X\|^2, \qquad \Tr (R_X)^2 = H \, \|X\|^4,
\end{equation}
with some constants $C,\, H$ (see, for instance, \cite[Ch.6.C]{B}). In particular, any
harmonic manifold is Einstein.

Flat and rank-one symmetric spaces are harmonic, as the isometry group of each of them
acts transitively on its unit tangent bundle. In \cite{Li}, Lichnerowicz conjectured that the
converse is true: any harmonic space is two-point homogeneous.

The Lichnerowicz Conjecture is known to be true for:
\begin{enumerate}[(a)]
    \item manifolds of positive scalar curvature (compact manifolds with finite fundamental group)
    and Ricci-flat manifolds \cite{S1, S2};
    \item compact manifolds of negative curvature (\cite{BCG}, using the result of \cite{FL} on asymptotic harmonicity);
    \item manifolds of dimension $n \le 5$ (\cite{Li,W} for $n \le 4$ and \cite{N} for $n=5$);
    \item symmetric spaces \cite{L},
\end{enumerate}
and also for some other classes of manifolds.

However, in 1992, Damek and Ricci \cite{DR} discovered an infinite series of harmonic
homogeneous non-compact spaces, which are, in general, non-symmetric, thus
disproving the Lichnerowicz Conjecture. A \emph{Damek-Ricci} space is a solvmanifold
(a solvable Lie group with a left invariant metric) whose Lie algebra is
constructed as follows. Let $\n = \z \oplus \mathfrak{u}$ be an orthogonal decomposition of a
nilpotent Lie algebra $\n$ with the center $\z$, and with
$[\mathfrak{u}, \mathfrak{u}] \subset \z$ ($\n$ is called two-step nilpotent).
For $Z \in \z$, define a skew-symmetric operator $J_Z: \mathfrak{u} \to \mathfrak{u}$ by
$\<J_Z X, Y\> = \<Z, [X, Y]\>$, for $X, Y \in \mathfrak{u}$.
A Lie algebra of a Damek-\-Ricci space is $\g = \n \oplus \mathfrak{a}$, with $\mathfrak{a}$ a
one-dimensional space orthogonal to $\n$, and with
$$
\ad_{A|\mathfrak{u}} = \tfrac{\la}{2}\, \id_{\mathfrak{u}},\quad \ad_{A|\z} =  \la \, \id_{\z},
\quad J_Z^2 = -\la \|Z\|^2 \id_{\mathfrak{u}},
$$
for all $Z \in \z$, where $A$ is a unit vector in $\mathfrak{a}$ and $\la$ a positive
constant (the operators $J_Z$ are constructed using Clifford modules).

Note that rank-one noncompact symmetric spaces (including the real hyperbolic space, if one allows
$\mathfrak{u} =0$) are specific cases of Damek-Ricci spaces. Recently, in \cite{BPR} it was proved
that a rank one three-step harmonic solvmanifold is a Damek-Ricci space (this strengthens the result
of \cite{D}). In \cite{H1}, it is proved that a homogeneous Hadamard manifold is harmonic if and
only if $\dim \mathfrak a = 1$ and the geodesic symmetries are volume preserving:
$V_{x, X}(t)=V_{x, -X}(t)$ for all $X \in T_xM^n$.

For further results on harmonic spaces and Damek-Ricci spaces we refer to
\cite{BTV, S2}.

\medskip

All the known harmonic spaces are homogeneous. This raises two questions:

1. Is a harmonic space necessarily homogeneous?

2. What are homogeneous harmonic spaces?

The first question was asked in \cite[Ch.4.5]{BTV}, and to the best of our knowledge,
is still open (for spaces of negative scalar curvature).
In this paper, we deal with the second one. Replacing the assumption of negativity
of the Ricci curvature by a stronger one, the nonpositivity of the sectional curvature,
we prove the following:

\begin{theorem*} A harmonic homogeneous manifold of nonpositive curvature is either flat,
or is locally isometric to a Damek-Ricci space.
\end{theorem*}

In Section~\ref{s:plan} we give the plan of the proof of the Theorem, the proof itself is contained
in Section~\ref{s:proof}.

\section{Plan of the proof}\label{s:plan}


By the result of \cite{AK}, a Ricci-flat homogeneous space is flat. In what follows
we therefore assume that the scalar curvature is negative.

As it follows from \cite{A, He} (see also \cite{AW1, AW2}), a homogeneous
manifold of non-positive curvature is a solvmanifold, that is, a solvable Lie
group $G$ with a left invariant metric (from this point on we denote the manifold $G$
instead of $M$). What is more, the orthogonal complement
$\mathfrak a$ to the nilradical $\n=[\g,\g]$ is an abelian subalgebra of the Lie algebra
$\g$ of $G$. By the result of \cite[Theorem 4.7]{H1}, $G$ is harmonic, only if $\dim \mathfrak a = 1$.

Furthermore, any harmonic manifold is Einstein, so by \cite{AW1} (see also \cite{H2}), $G$ is
isometric to a solvmanifold of \emph{Iwasawa type}, which, when $\dim \mathfrak a = 1$,
means that a unit vector $A \in \mathfrak a$ can be chosen in such a way that the operator
$D:=\ad_{A|\n}$ is symmetric and positive definite.

So, what we really have to prove, is that a harmonic solvmanifold
\begin{enumerate}[(i)]
    \item of Iwasawa type,
    \item of rank $1$ (i.e., $\dim \mathfrak a = 1$), and
    \item of non-positive curvature
\end{enumerate}
is either flat, or is a Damek-Ricci space. Note that we already used the assumption of
nonpositivity of the curvature, but it will still be needed once again further in the proof.

We have an orthogonal decomposition of $\n$ on the eigenspaces of $D$. Namely, for $\a >0$, let
$\n_\a = \{X \in \n: DX = \a X\}$, so that
$D_{|\n_\a} = \a \; \id_{|\n_\a}$, and let $\Delta =\{\a: \dim \n_\a >0\}$.

As $D$ is a derivation of $\n, \; [\n_\a, \n_\b] \subset \n_{\a+\b}$. In particular, if $\la = \max \Delta$,
the biggest eigenvalue of $D$, then the eigenspace $\n_\la$ lies in the center $\z$ of $\n$.

We start with choosing (and fixing) a unit vector $Z \in \n_\la \subset \z$ and considering a geodesic
$\G_\phi$ of $G$ passing through $e \in G$ in the direction $A \cos \phi + Z \sin \phi$. When $\phi=0$,
the geodesic $\G_0$ (\emph{the abelian geodesic}) is a one-dimensional subgroup of $G:\; \G_0(t) = \exp_e(tA)$.
This is no longer true for an arbitrary $\phi$. However, the fact that the subalgebra $\Span (A, Z)$ is tangent
to a hyperbolic plane, which is totally geodesic in $G$,
makes it possible to find the equation of the geodesic $\G_\phi$ explicitly, for any $\phi$, and significantly
simplifies the equation for Jacobi fields along $\G_\phi$ (Lemmas~\ref{l:geodesic}, \ref{l:Jacobi}, \ref{l:Taylor}).

As the solvmanifold $G$ is harmonic, the volume density function $V(t, \phi):=V_{e,A \cos \phi + Z \sin \phi}(t)$
along $\G_\phi$ must not depend of $\phi$. Moreover, along the abelian geodesic
$\G_0, \; V(t, 0) = \prod_{\a \in \Delta} (\a^{-1} \sinh \a t)^{\dim \n_\a}$.

Considering the Taylor expansion of $V(t, \phi)$ at $\phi = 0$ we find that
$\tfrac {d}{d\phi} (V(t, \phi))_{\phi=0} \equiv 0$ (in fact, $V(t, \phi)$ is an even function of $\phi$),
but the condition $\tfrac {d^2}{d\phi^2} (V(t, \phi))_{\phi=0} = 0$
gives nontrivial restrictions on the eigenvalue structure (Lemma~\ref{l:rootspace}):
$$
\a \dim \n_\a = (\la-\a) \dim \n_{\la-\a} \quad (= (2 \la)^{-1} \Tr (J_Z^{} J_Z^t)_{|\n_\a}),
\quad \text{for $\a \in \Delta \setminus \{\la\}$},
$$
where $J_ZX = \ad^*_XZ$ for $X \in \n$. It follows that the eigenvalues other
than $\la$ and $\tfrac12 \la$ (if the latter is an eigenvalue) come in pairs: $\a, \; \la - \a$. What
is more, the above equation, together with the fact that the curvature is nonpositive, implies
that $\Delta \subset [\tfrac13 \la, \tfrac23 \la] \cup \{\la\}$.

Using \cite[Theorem 4.14]{H2} we further narrow the set $\Delta$ in Lemma~\ref{l:rationality}:
$\Delta \subset \{\tfrac13 \la, \tfrac12 \la, \tfrac23 \la, \la\}$.

Next, in Lemma~\ref{l:la/2}, we apply the second Ledger formula to prove that
$J_Z^2X = -\la^2 \|Z\|^2 X$ for all $X \in \n_{\la/2}, \; Z \in \n_\la$, and that the cases
$\tfrac12 \la \in \Delta$ and $\{\tfrac13 \la, \tfrac23 \la\} \cap \Delta \ne \varnothing$ are
mutually exclusive.
Then either $\Delta \subset \{\tfrac12 \la, \la\}$, which leads to
Damek-Ricci spaces, or $\Delta = \{\tfrac13 \la, \tfrac23 \la, \la\}$.
To finish the proof it remains to show
that in the latter case, the volume density function $V(t, \phi)$ of such, quite a specific, solvmanifold
depends on $\phi$. This is done in Lemma~\ref{l:notla/3}: the coefficient of $t^9$ in the Taylor expansion
of $V(t, \phi)$ appears to be nonconstant.

\section{Proof}\label{s:proof}

\subsection{Notations and basic facts}\label{ss:facts}


Let $G$ be a solvmanifold of Iwasawa type, $\dim G = n$, and $\g = \mathfrak a \oplus \n$ be an orthogonal decomposition of its Lie
algebra, with the nilradical $\n = [\g, \g]$ and $\dim \mathfrak a =1$.
The operator $D = \ad_{A|\n}$ is symmetric and positive definite (for one of two possible choices of a unit vector
$A \in \mathfrak a$).

For $\a >0$, let $\n_\a = \{X \in \n: DX = \a X\}$. Denote $\Delta =\{\a: \dim \n_\a >0\}$, and let $\la = \max \Delta$.
We have $[\n_\a, \n_\b] \subset \n_{\a+\b}$, so in particular, $\n_\la \subset \z$, the center of $\n$.

For $\a \in \Delta$, let $n_\a = \dim \n_\a$, the multiplicity of the eigenvalue $\a$ of $D$. We denote
$\la_1, \la_2, \ldots, \la_{n-1}$ the eigenvalues of $D$ counting multiplicities, so that
$(\la_1, \la_2, \ldots, \la_{n-1})$ contains $n_\a$ copies of each $\a \in \Delta$.
A particular labelling of the $\la_i$'s will be specified further in the proof.

For $Z \in \n_\la$, define a skew-symmetric operator $J_Z: \g \to \g$ by $\<J_Z U, V\> = \<Z, [U,V]\>$.

For left invariant vector fields, the connection and the curvature are given by
\begin{equation}
\nabla_VW= U(V,W)+ \tfrac12 [V, W], \quad \text{where} \quad \<U(V,W),Y\> = \tfrac12(\<V, [Y, W]\> + \<W, [Y, V]\>), \label{eq:nabla}
\end{equation}
\begin{multline}\label{eq:R}
R(X, Y, Y, X) = \|U(X, Y)\|^2 - \<U(X, X), U(Y, Y)\> -\tfrac34 \|[X, Y]\|^2 \\
    -\tfrac12 \<[X, [X, Y]], Y\> -\tfrac12 \<[Y, [Y, X]], X\>,
\end{multline}
respectively. As usual, we identify left invariant vector fields on $G$ with their values at $e$.

In the following Lemma, we collect some simple facts to be used further in the proof.

\begin{lemma}\label{l:facts} Let $Z \in \n_\la \subset \z$ be a unit vector and $\tn = \n \cap Z^\perp$.
Then

\emph{1.} $J_Z Z= -\la A, \; J_Z A = \la Z, \; J_Z \tn \subset \tn, \; J_Z \n_\a \subset \n_{\la-\a}$. In particular,
$J_Z \n_\a = 0$, when $\la-\a \notin \Delta$.

\emph{2.} The restrictions of the symmetric operators $J_Z^2$ and $D$ to their invariant subspace $\tn$ commute.
For every $\a \in \Delta$, there is an orthogonal decomposition
$\n_\a = (\n_\a \cap \Ker J_Z) \oplus (\n_\a \cap (\Ker J_Z)^\perp)$.

\emph{3.} $\nabla_A U = 0$ for any left invariant vector field $U$ \emph{(}in particular, $\exp (tA)$ is a geodesic\emph{)}.

\emph{4.} The distribution $\mathcal{D}= \Span(A,Z)$ on $G$ is totally geodesic; its integral manifolds are isometric
to the hyperbolic plane of curvature $-\la^2$.
\end{lemma}

\begin{proof} 1. The first two equations follow from the definition of $J_Z$. As $\Span(A,Z)$ is an invariant
subspace of $J_Z$, its orthogonal complement $\tn$ also is. Moreover, for
$X \in \n_\a, Y \in \n_\b,\; \<J_Z X, Y\> = \<Z, [X,Y]\>$, which can be nonzero only when $[X,Y] \in \n_\la$,
that is, when $\a + \b =\la$.

2. 
The first claim follows from the fact that $J_Z \n_\a \subset \n_{\la-\a}$ and $D_{|\n_\a} = \id_{|\n_\a}$.
To prove the second one, notice that $J_Z$ and $J_Z^2$ have the same kernel and cokernel.

3. Follows from \eqref{eq:nabla}.

4. The distribution $\mathcal{D}$ is integrable, as $[A, Z] = \la Z \in \mathcal{D}$.
Its integral submanifolds are totally geo\-desic, since
$\nabla_A A = \nabla_A Z = 0,\; \nabla_Z Z = \la A, \; \nabla_Z A = - \la Z$. By \eqref{eq:R},
$R(A, Z, Z, A) = - \la^2$.

\end{proof}

In the next Lemma we consider a geodesic $\G_\phi(t)$ of $G$ lying in the hyperbolic plane introduced in
assertion~4 of Lemma~\ref{l:facts} (dot stands for $\tfrac{d}{dt}$).

\begin{lemma}\label{l:geodesic}
Let $\G_\phi(t)$ be a unit speed geodesic of the solvmanifold $G$, which passes
through $e$ in the direction $A \cos \phi + Z \sin \phi$, where $Z \in \n_\la$ is a unit vector. Then

\emph{1.} The velocity vector $\dot \G_\phi(t)$ of $\G_\phi(t)$ is given by
    \begin{gather}
    \dot \G_\phi(t)= q A + \Phi Z, \quad \text{where} \notag \\
    q = \frac{-\sinh \la t + \cos \phi \cosh \la t}{\cosh \la t - \cos \phi \sinh \la t}, \quad
    \Phi = \frac{\sin \phi}{\cosh \la t - \cos \phi \sinh \la t}, \quad q^2 + \Phi^2 = 1. \label{eq:qandPhi}
    \end{gather}

\emph{2.} The vector bundle $\tn$ along $\G_{\phi}$ is orthogonal to $\dot \G_\phi(t)$ and is invariant with
    respect to $\nabla_{\dot \G_\phi(t)}$. For a left invariant vector field $X \in \tn$,
    \begin{equation}\label{eq:nablaX}
    \nabla_{\dot \G_\phi(t)}X = -\tfrac12 \Phi \, J_Z X,
    \end{equation}

\emph{3.} The bundle $\tilde \n$ is invariant with respect to the Jacobi operator $R_{\dot \G_\phi(t)}$, and
    \begin{equation}\label{eq:RdotG}
    R_{\dot \G_\phi(t)| \tn } = (- q^2 D^2 -\tfrac14 \Phi^2J_Z^2 - \la \Phi^2 D + \tfrac12 q \Phi \, [D,J_Z])_{| \tn }.
    \end{equation}

\end{lemma}

\begin{proof}1. As it follows from assertion 4 of Lemma~\ref{l:facts}, the geodesic $\G_{\phi}$ lies
in a totally geodesic hyperbolic plane $\mathbb{H}(-\la^2) = \exp (\Span(A, Z)) \subset G$ passing through $e$.
It follows that $\dot \G_\phi(t)= q A + \Phi Z$, with some functions $q=q(\phi, t)$ and $\Phi = \Phi(\phi, t)$
satisfying $q^2 + \Phi^2 =1$. The fact that $\nabla_{\dot \G_\phi(t)}\dot \G_\phi(t)=0$ implies that
$\dot \Phi = \la \Phi q,\; \dot q + \la \Phi^2 = 0$. Solving this, subject to the initial conditions
$q(\phi, 0) = \cos \phi,\; \Phi (\phi, 0) = \sin \phi$, we get what required.

It is not difficult to find the equation of the geodesic $\G_\phi(t)$ explicitly, but we
won't need it.

2. This follows immediately from assertion 1 and the fact that $\nabla_AX=0$ (assertion 3 of Lemma~\ref{l:facts}) and
$\nabla_ZX = -\tfrac12J_ZX$ (equation \eqref{eq:nabla}).

3. The fact that $\tn$ is an invariant subspace of $R_{\dot \G_\phi(t)}$ follows from assertion 4
of Lemma~\ref{l:facts}. Equation \eqref{eq:RdotG} is a direct consequence of \eqref{eq:R}.
\end{proof}

\subsection{Jacobi fields and the volume density function}\label{ss:jacobi}

In this subsection, we find the equation of the Jacobi fields along the geodesic $\G_\phi(t)$, and then calculate
$\tfrac{d^2}{d\phi^2} V(t, \phi)$ at $\phi = 0$.

By assertion 2 of Lemma~\ref{l:facts}, the symmetric operators $D_{| \tn }$ and $J^2_{Z| \tn}$ commute,
so we can choose an orthonormal basis
$X_1, \ldots, X_k, Y_1, Y_2, \ldots, Y_{2p-1}, Y_{2p}, V_1, V_2, \ldots, V_{2m-1}, V_{2m} \; (k + 2p + 2m = n-2)$
in $\tn$ in such a way that
\begin{align}\label{eq:basisX}
D(X_j) &= \la_j X_j,& J_Z X_j &= 0, && j = 1, \ldots, k,\\ \label{eq:basisY}
D(Y_{2i-1}) &= \tfrac12 \la Y_{2i-1}, \; D(Y_{2i}) = \tfrac12 \la Y_{2i}, & J_Z Y_{2i-1} &= a_i Y_{2i}, \; J_Z Y_{2i} = -a_i Y_{2i-1},
&& i = 1, \ldots, p,\\
D(V_{2l-1}) &= \la_l V_{2l-1}, \; D(V_{2l}) = (\la -\la_l) V_{2l}, & J_Z V_{2l-1} &= b_l V_{2l}, \; J_Z V_{2l} = -b_l V_{2l-1},
&& l = 1, \ldots, m, \label{eq:basisV}
\end{align}
with $a_i, b_l \ne 0$ and $\la_l < \tfrac12 \la$ (note that the $\la$'s here need not to be distinct: for instance, some
of the $\la_j$'s, with $j=1, \ldots, k$ can be equal to $\tfrac12 \la$, or to one of the $\la_l$ or $\la-\la_l, \; l = 1, \ldots, m$).

The vector bundle $\tn$ along $\G_\phi$ splits into orthogonal sum of one-dimensional subbundles
$\Span(X_j)$ and two-dimensional subbundles
$\Span(Y_{2i-1}, Y_{2i}),\; \Span(V_{2l-1}, V_{2l})$, each of
which is invariant with respect to both $\nabla_{\dot \G_\phi(t)}$ and the Jacobi operator $R_{\dot \G_\phi(t)}$
(assertions 2 and 3 of Lemma~\ref{l:geodesic}). So the Jacobi equation along $\G_\phi$ splits into a set
of single equations and pairs of equations. We have the following:

\begin{lemma}\label{l:Jacobi}
The volume density function $V(t,\phi)$ along the geodesic $\G_\phi$ has the form
\begin{equation}\label{eq:volume}
    V(t, \phi) = \frac{\sinh \la t}{\la} \prod_{j=1}^k x_j(t, \phi) \prod_{i=1}^p y^2_i(t, \phi) \prod_{l=1}^m \det v_l(t,\phi),
\end{equation}
where the functions $x_j(t, \phi)$ and $y_i(t, \phi)$ are determined by
\begin{align}\label{eq:x}
    \ddot x_j &= (\la_j^2 + (\la - \la_j) \la_j\Phi^2) x_j, & x_j (0, \phi) &= 0, & \dot x_j (0, \phi) &= 1,\\
    \ddot y_i &= \bigl(\tfrac14 \la^2 + \tfrac{\la^2 - a_i^2}4 \Phi^2\bigr) y_i, & y_i (0, \phi) &= 0, & \dot y_i (0, \phi) &= 1,
    \label{eq:y}
\end{align}
and the $2 \times 2$-matrix $v_l(t,\phi)$ satisfies
\begin{equation}\label{eq:v}
    \ddot v_l + b_l \Phi \left(\begin{array}{cc} 0 & 1 \\ -1 & 0 \\ \end{array}\right) \dot v_l +
    \left(\begin{array}{cc} -(q^2 \la_l^2 + \la\la_l \Phi^2) & q \Phi \; b_l \la_l' \\
    - q \Phi \; b_l \la_l & -(q^2 \la_l'^2 + \la\la_l' \Phi^2) \\ \end{array}\right) v_l = 0, \quad
    \begin{array}{l} v_l (0, \phi) = 0, \\ \dot v_l (0, \phi) = I_2, \\ \end{array}
\end{equation}
where $\la_l' = \la-\la_l$, and $I_2$ is the $2 \times 2$ identity matrix. If $\phi = 0$, then $\Phi \equiv 0,\; q \equiv 1$, and
so $x_j (t, 0) = \tfrac 1{\la_j} \sinh \la_j t,\;  y_i (t, 0) = \tfrac 2{\la} \sinh \tfrac{\la}{2} t,\;
v_l (t, 0) = \diag(\tfrac{1}{\la_l} \sinh \la_l t, \tfrac 1{\la_l'} \sinh \la_l' t)$
and
\begin{equation}\label{eq:volume0}
    V(t, 0) = \frac{\sinh \la t}{\la}  \prod_{j=1}^k \frac{\sinh \la_j t}{\la_j}  \prod_{i=1}^p \Bigl(\frac {\sinh \tfrac{\la}{2}t}{\tfrac{\la}{2}}
    \Bigr)^2
    \prod_{l=1}^m \frac {\sinh \la_l t \sinh \la_l' t}{\la_l \la_l'}  = \prod_{\a \in \Delta} \left(\frac{\sinh \a t}{\a} \right)^{n_\a}.
\end{equation}
\end{lemma}

\begin{proof} Consider Jacobi fields along the geodesic $\G_\phi$, orthogonal to $\dot \G_\phi(t)$ and vanishing at $t=0$.
First of all, as
the distribution $\mathcal{D}=\Span (A, Z)$ is tangent to the totally geodesic hyperbolic plane of curvature $-\la^2$
(assertion 4 of Lemma~\ref{l:facts}), there is a Jacobi field along $\G_\phi(t)$ lying in $\mathcal{D}$, whose contribution
to the volume density function is $\frac{\sinh \la t}{\la}$. Moreover, all the Jacobi fields whose derivative at $t=0$ is
orthogonal to $\mathcal{D}$, remain orthogonal to $\mathcal{D}$ for all $t$.
From \eqref{eq:RdotG} we have:
\begin{equation}\label{eq:RdotGX}
    R_{\dot \G_\phi(t)} X =
\begin{cases}
    (-\a^2- \a\,(\la-\a) \Phi^2)X, & X \in \n_\a \cap \Ker J_Z, \, X \perp Z; \\
    -\tfrac14 (1 + \Phi^2) \la^2 X - \tfrac14 \Phi^2 J_Z^2X, & X \in \n_{\la/2} \cap (\Ker J_Z)^\perp;\\
    (-q^2\a^2 - \la \a \Phi^2) X - \tfrac14 \Phi^2 J_Z^2X + \tfrac12 q \Phi (\la-2\a)J_ZX, &
    X \in \n_\a \cap (\Ker J_Z)^\perp, \; \a \ne \tfrac{\la}{2}.
\end{cases}
\end{equation}

Let now $X_j$ be a unit left invariant vector field defined by \eqref{eq:basisX}. As $\Span(X_j)$ is invariant with respect to
both $R_{\dot \G_\phi(t)}$ and $\nabla_{\dot \G_\phi(t)}$, there is a Jacobi field scalar proportional to $X_j$. From
\eqref{eq:nablaX} and the first line of \eqref{eq:RdotGX} we find that the scale factor
$x_j(t, \phi)$ satisfies \eqref{eq:x}.

Similarly, the left invariant vector bundle $\Span(V_{2l-1}, V_{2l})$ defined by \eqref{eq:basisV} is invariant with respect to
both $R_{\dot \G_\phi(t)}$ and $\nabla_{\dot \G_\phi(t)}$, hence it contains two linearly independent Jacobi fields vanishing
at $t = 0$. From \eqref{eq:nablaX}, the fact that $\dot \Phi = \la q \Phi$ (assertion 1 of Lemma~\ref{l:geodesic}) and the third
line of \eqref{eq:RdotGX} it follows that the $2 \times 2$-matrix $v_l(t,\phi)$ of their components with respect to the orthonormal
basis $V_{2l-1}, V_{2l}$ satisfies \eqref{eq:v}.

The same arguments for the left invariant vector bundle $\Span(Y_{2i-1}, Y_{2i})$ defined by \eqref{eq:basisY}
give two linearly independent Jacobi fields, whose coordinates in the orthonormal basis $Y_{2i-1}, Y_{2i}$ satisfy
the equation
\begin{equation*}
    \ddot w_i + a_i \Phi J \dot w_i + \tfrac12 a_i \dot \Phi J w_i - \tfrac14 \la^2 (1 + \Phi^2) w_i = 0,
    \quad \text{with $w_i (0, \phi) = 0, \quad \dot w_i (0, \phi) = I_2$},
\end{equation*}
where $J=\left(\begin{smallmatrix} 0 & 1 \\ -1 & 0 \\ \end{smallmatrix}\right)$.
Substituting $w_i= \exp(-\tfrac12 a_i \int_0^t \Phi J dt) z_i$ we get
$\ddot z_i =\tfrac14 (\la^2 + (\la^2 - a_i^2) \Phi^2) z_i$, with the initial conditions
$z_i (0, \phi) = 0, \; \dot z_i (0, \phi) = I_2$. It follows that the $2 \times 2$-matrix $z_i$ is proportional to
the identity matrix, and $\det z_i = y_i^2$, with the function $y_i$ given by \eqref{eq:y}.
\end{proof}

\begin{lemma}\label{l:Taylor}
Let $x_j(t, \phi),\; y_i(t, \phi)$ and $v_l(t, \phi)$ satisfy \eqref{eq:x}, \eqref{eq:y} and \eqref{eq:v}, respectively. Then, in
a neighborhood of $(t, \phi) = (0, 0)$,
\begin{gather*}
    x_j(t, \phi) = \frac{\sinh \la_j t}{\la_j} \bigl(1 + \frac{\phi^2}2 e^{\la t} \sinh \la t \;\hat x_j(t) + o(\phi^2)\bigr), \;
    y_i(t, \phi) = \frac{\sinh \tfrac{\la}{2} t}{\la/2} \bigl(1 + \frac{\phi^2}2 e^{\la t} \sinh \la t \;\hat y_i(t) + o(\phi^2)\bigr)\\
    \det v_l(t, \phi) = \frac {\sinh \la_l t}{\la_l} \; \frac {\sinh \la_l' t}{\la_l'}
    \bigl(1 + \frac{\phi^2}2 e^{\la t} \sinh \la t \; \hat v_l(t) + o(\phi^2)\bigr),
\end{gather*}
where
\begin{align}\label{eq:secondderx}
    \hat x_j(t) &= \frac{\la_j}{\la(\la+\la_j)}(\la \coth\la t -\la_j \coth \la_j t),\\ \label{eq:seconddery}
    \hat y_i(t) &= \frac{\la^2 - a_i^2}{6 \la^2}(2 \coth\la t - \coth \tfrac{\la}{2}t),\\ \label{eq:secondderv}
    \hat v_l(t) &= \frac{2(\la^2 + 2\la_l\la_l')-3 b_l^2}{2 (\la+\la_l)(\la+\la_l')} \coth\la t
    + \frac{\la_l(b_l^2 - 2\la \la_l)}{2 \la^2 (\la+\la_l)} \coth \la_l t
    + \frac{\la_l'(b_l^2 - 2\la \la_l')}{2 \la^2 (\la+\la_l')} \coth \la_l' t,
\end{align}
and so
\begin{equation}\label{eq:volumeTaylor}
    V(t, \phi) = V(t, 0) \Bigl(1 + \frac{\phi^2}2 e^{\la t} \sinh \la t \;
    \Bigl(\sum_{j=1}^k \hat x_j(t) + 2 \sum_{i=1}^p \hat y_i(t) + \sum_{l=1}^m \hat v_l(t)\Bigr) +o(\phi^2)\Bigr).
\end{equation}
\end{lemma}

\begin{proof} Both the equation \eqref{eq:x} and \eqref{eq:y} have the form
$$
\ddot f= (c_1^2 + c_2 \Phi^2) f, \quad f(0, \phi) = 0, \quad \dot f (0, \phi) = 1,
$$
with some constants $c_1 \in (0, \la)$ and $c_2$. Expanding $f(t, \phi) = f_0(t) + \phi f_1(t) + \tfrac{\phi^2}{2} f_2(t) + o(\phi^2)$,
using the fact that $\Phi^2(t,\phi)= \phi^2 e^{2\la t} + o(\phi^2)$ (from \eqref{eq:qandPhi}), and substituting to the equation above, we
find:
$$
f_0(t) = \frac{\sinh c_1 t}{c_1}, \quad f_1(t)=0, \quad
f_2(t)= \frac{c_2 e^{\la t}}{\la c_1 (\la^2 - c_1^2)} (\la \sinh c_1 t \cosh \la t - c_1 \cosh c_1 t \sinh \la t),
$$
which gives the required expansions for $x_j(t, \phi)$ and $y_i(t, \phi)$.

Similarly, for the matrix $v_l(t, \phi)$, the solution to \eqref{eq:v}, we have:
\begin{multline*}
v_l(t, \phi) = \left(%
\begin{array}{cc}
  \!\!\tfrac{\sinh \la_l t}{\la_l} & 0 \\
  0 & \!\!\tfrac{\sinh \la_l' t}{\la_l'}\!\! \\
\end{array}%
\right) + \, \phi \, \frac{b_l}{2 \la \la_l^{} \la_l'}
\left(%
\begin{array}{cc}
  0 & \!\!(\la_l' \sinh \la_l^{} t - e^{\la t} \la_l^{} \sinh \la_l' t)\\
  \!\!(\la_l^{} \sinh \la_l' t - e^{\la t} \la_l' \sinh \la_l^{} t)\!\! & 0 \\
\end{array}%
\right)  \\ - \frac{\phi^2}{2} \frac{b_l^2 - \la_l^{}\la_l'}{\la \la_l^{}\la_l'} \, e^{\la t} \,
\left(%
\begin{array}{cc}
  \frac{\la \cosh \la t \sinh \la_l t - \la_l  \sinh \la t \cosh \la_l t}{\la + \la_l^{}}  & 0 \\
  0 & \frac{\la \cosh \la t \sinh \la_l' t - \la_l'  \sinh \la t \cosh \la_l' t}{\la + \la_l'}  \\
\end{array}%
\right)
+ o(\phi^2),
\end{multline*}
and the expansion for $\det v_l(t, \phi)$ follows.
\end{proof}
Note (although we won't use this) that the equation \eqref{eq:x} can be solved in the closed form by the reduction of order:
$x(t, \phi) = (\cosh \la t - \cos \phi \sinh \la t)^{\la_j/\la} \int_0^t (\cosh \la u - \cos \phi \sinh \la u)^{-2\la_j/\la}du$.

\subsection{The eigenvalue set}\label{ss:root}

In this subsection, we show that the expansion of the volume density function obtained in Lemma~\ref{l:Taylor},
gives strong restrictions on the eigenvalue set $\Delta$. After few steps, this, together with the nonpositivity of
the curvature, leaves only two possibilities: either $\Delta = \{\la\}$, or
$\Delta = \{\tfrac12\la, \la\}$, with $J_Z^2X = -\la^2 \|Z\|^2 X$ for all $X \in \n_{\la/2}, \; Z \in \n_\la$.
This completes the proof of the Theorem as the former one corresponds to a hyperbolic space and the latter one to
a Damek-Ricci space.

\begin{lemma}\label{l:rootspace}
\emph{1.} We have:
\begin{gather}\label{eq:ana}
\a \, n_\a = (\la-\a) \, n_{\la-\a} = \frac{1}{2 \la} \sum_{l: \la_l = \a} b_l^2 = \frac{1}{2 \la} \Tr (J_Z^{}J_Z^t)_{|\n_\a},
\quad \text{for any $\a \in \Delta, \; \a \ne \la, \la/2$} \\
\la^2 n_{\la/2} = 2 \sum_{i=1}^p a_i^2 =  \Tr (J_Z^{}J_Z^t)_{|\n_{\la/2}},\quad \text{if $\la/2 \in \Delta$.} \label{eq:la/2}
\end{gather}
In particular, for any $\a \in (0, \la), \quad \a \in \Delta \Leftrightarrow \la - \a \in \Delta$;
for $\a \in \Delta \setminus \{\la\}, \quad \n_\a \cap (\Ker J_Z)^\perp \ne 0$.

\emph{2.} $\Delta \subset [\tfrac13 \la, \tfrac23 \la] \cup \{\la \}$.

\emph{3.} 
$(J_Z^{}J_Z^t)_{|\n_{2\la/3}} = \tfrac43 \la^2 \id_{|\n_{2\la/3}}$.
\end{lemma}
\begin{proof}
1. The harmonicity of $G$ implies that $V(t, \phi) = V(t, 0)$, for all $\phi$, and so, by \eqref{eq:volumeTaylor},
\begin{equation}\label{eq:sum0}
\sum_{j=1}^k \hat x_j(t) + 2 \sum_{i=1}^p \hat y_i(t) + \sum_{l=1}^m \hat v_l(t)=0,
\end{equation}
with $\hat x_j(t),\; \hat y_i(t)$ and $\hat v_l(t)$
given by \eqref{eq:secondderx}, \eqref{eq:seconddery} and \eqref{eq:secondderv}, respectively.
The expression on the left hand side is a linear combination of the $\coth$'s. Note that
for any finite set of positive numbers $\mu_1 > \mu_2 > \ldots > \mu_N$, the functions $\coth (\mu_s t)$
are linearly independent. For if some linear combination $\sum_s c_s\coth (\mu_s t)$ vanishes
for all $t \in \mathbb{R} \setminus \{ 0 \}$, it will still vanish for all $t \in \mathbb{C}$,
outside the union of the poles of the meromorphic functions $\coth (\mu_s t)$. However, the
function $\coth (\mu_1 t)$ goes to infinity when $t \to \pi i/\mu_1$, while all the others stay
bounded. So $c_1=0$, and, inductively, all the $c_s$'s are zeros.

Equating the coefficients of $\coth \a t$ in \eqref{eq:sum0} to zero, we find:

For $\a \in \Delta, \; \a <\tfrac{\la}{2}: \;
\sum_{j: \la_j = \a} \frac{-\a^2}{\la(\la+\a)} + \sum_{l: \la_l = \a} \frac{\a (b_l^2 - 2\la \a)}{2 \la^2 (\la+\a)}=0$,
and so $2\a\la (\#\{j: \la_j = \a\} + \#\{l: \la_l = \a\}) = \sum_{l: \la_l = \a} b_l^2$. The expression in the brackets
is $\dim (\n_\a \cap \Ker J_Z) + \dim (\n_\a \cap (\Ker J_Z)^\perp) = n_\a$ (by assertion 2 of Lemma~\ref{l:facts}). So
$$
2\a \la n_\a = \sum_{l: \la_l = \a} b_l^2 = \Tr (J_Z^{}J_Z^t)_{|\n_\a},
$$
the latter equation follows from the definition of the $b_l$'s. The same formula remains true also for
$\a \in \Delta \cap (\tfrac{\la}{2}, \la)$ (replace $\la_l$ by $\la_l'$ in the above computation).
Since $J_Z\n_\a \subset \n_{\la-\a},\; J_Z\n_{\la-\a} \subset \n_\a$, and the operator $J_Z$ is skew symmetric,
$\Tr (J_Z^{}J_Z^t)_{|\n_\a} = \Tr (J_Z^{}J_Z^t)_{|\n_{\la-\a}}$. We get
$$
2\a \la n_\a = 2(\la-\a) \la n_{\la-\a} = \Tr (J_Z^{}J_Z^t)_{|\n_\a},
$$
which proves \eqref{eq:ana}.
Next, for $\a = \la/2$, the coefficient of $\coth \tfrac{\la}2 t$ on the left hand side of \eqref{eq:sum0} is
$-\tfrac16\#\{j: \la_j = \la/2\} - 2 \sum_{i=1}^p \frac{\la^2 - a_i^2}{6 \la^2}=0$, so
$$
\la^2 n_{\la/2} = 2 \sum_{i=1}^p a_i^2 =  \Tr (J_Z^{}J_Z^t)_{|\n_{\la/2}},
$$
which is \eqref{eq:la/2}. It can be seen that the coefficient of $\coth \la t$ does not give anything new.

The fact that $\n_\a \cap (\Ker J_Z)^\perp \ne 0$ for $\a \ne \la$ follows from (\ref{eq:ana}, \ref{eq:la/2}),
as $a_i, b_l \ne 0 $.

2. For $\a \in \Delta,\; \a \ne \la$, let $i_\a = \dim (\n_\a \cap (\Ker J_Z)^\perp)$. As it follows
from assertion 1, $i_\a > 0$, and obviously  $i_\a \le n_\a$. Moreover, as $J_Z$ is skew-symmetric
and as $J_Z \n_\a \subset \n_{\la-\a}$,  we have $i_\a = i_{\la-\a}$, so
$$
0 < i_\a = i_{\la-\a} \le n_\a, n_{\la-\a}.
$$
For a unit vector $V \in \n_\a \cap (\Ker J_Z)^\perp, \quad R(Z, V, V, Z) = \tfrac14 \|J_ZV\|^2 - \la \, \a$
(from the third line of \eqref{eq:RdotGX}, with $q=0, \Phi =1$).
As the curvature is nonpositive, $\|J_ZV\|^2 \le 4 \la \, \a$, and so, by \eqref{eq:ana},
$$
\a \, n_\a = (\la-\a) \, n_{\la-\a} = (2 \la)^{-1} \Tr (J_Z^{}J_Z^t)_{|\n_\a}
= (2 \la)^{-1} \Tr (J_Z^{}J_Z^t)_{|\n_\a \cap (\Ker J_Z)^\perp} \le 2 \, \a \, i_\a.
$$
Hence $(\la-\a)\, i_\a \le (\la-\a) n_{\la-\a} \le 2 \, \a \, i_\a$, which shows that $\a \ge \la/3$.
The fact that $\a \le 2 \la/3$ follows from assertion 1.

3. For $\a = \la/3$, the above inequalities become equations, and we obtain:
\begin{equation}\label{eq:ila}
\begin{split}
n_{\la/3} &= 2 \, i_{\la/3}, \quad (J_Z^2)_{|\n_{\la/3} \cap (\Ker J_Z)^\perp} = -\tfrac{4}{3} \la^2 \, \id_{\n_{\la/3} \cap (\Ker J_Z)^\perp}, \\
n_{2\la/3} &= i_{\la/3},\hphantom{2 \,} \quad (J_Z^2)_{|\n_{2\la/3}} = -\tfrac{4}{3} \la^2 \, \id_{\n_{2\la/3}}
\end{split}
\end{equation}
(the latter equation follows from the fact that $n_{2\la/3} = i_{\la/3}$, so $\n_{2\la/3} \perp \Ker J_Z$).
\end{proof}

\begin{lemma}\label{l:rationality}
$\Delta \subset \{\tfrac13 \la, \tfrac12 \la, \tfrac23 \la , \la \}$.
\end{lemma}
\begin{proof} We use equation (4.1) from the proof of Theorem 4.14 of \cite{H2}, in
slightly different notations.

Let $\Delta = \{\a_1, \ldots, \a_N\}$ and $c = \Tr D^2 /\Tr D$. In an orthonormal basis $\{e_1, \ldots, e_N\}$ for
$\mathbb{R}^N$, let $F$ be a set of vectors $\{e_i + e_j - e_k : \, \a_i + \a_j = \a_k, \; i \le j\}$.
Then
$$
(n_{\a_1}(c-\a_1), \ldots, n_{\a_N}(c-\a_N))^t \in \Span (F).
$$
Now assume there exists $\a \in \Delta$ such that $\a \in (\la/3, \la/2)$. Choose the labelling in such a way
that $\a_1 = \a,\; \a_2 = \la-\a, \; \a_N = \la$. Then (by assertion 2 of Lemma~\ref{l:rootspace}) the only relation of
the form $\a_i + \a_j = \a_k$ involving $\a$ is $\a_1 + \a_2 = \a_N$ (and the same is true for $\la -\a$). So $F$ contains a vector
$(1,1, 0, \ldots, 0, -1)^t$, and the first two coordinates of all the other vectors from $F$ are zeros. Hence
$n_\a(c-\a) = n_{\la-\a} (c-(\la-\a))$. This contradicts the fact that $\a \, n_\a= (\la-\a) \, n_{\la-\a}$,
which follows from \eqref{eq:ana}.
\end{proof}

\begin{lemma}\label{l:la/2}
\emph{1.} For any $X \in \n_{\la/2}, \; Z \in \n_\la$, we have $J_Z^2X = -\la^2 \|Z\|^2 X$.

\emph{2.} If $\la/2 \in \Delta$, then $\Delta = \{\la/2,\, \la\}$. 
\end{lemma}
\begin{proof} 1. Let $Z \in \n_\la$ be a unit vector.
From the second Ledger formula \eqref{eq:ledger},
$\Tr R_Z^2 = \Tr R_A^2$. As $R_A = - D^2$ (by \eqref{eq:RdotG}, with $\Phi = 0, \, q=1$),
\begin{equation}\label{eq:TrR2}
\Tr R_A^2 = \la^4 n_\la + \tfrac{\la^4 }{16} n_{\la/2} + \tfrac{\la^4}{81} n_{\la/3} + \tfrac{16 \la^4}{81} n_{2 \la/3} =
\la^4 n_\la + \tfrac{\la^4 }{16} n_{\la/2} + \tfrac{2 \la^4}{9} n_{2 \la/3},
\end{equation}
where we used the fact that $n_{\la/3} = 2 n_{2 \la/3}$.

On the other hand, from \eqref{eq:RdotGX}, with $q = 0, \, \Phi = 1$ and \eqref{eq:ila},
$$
R_ZX=
\begin{cases}
- \la^2 A & X = A, \\
- \la^2 X & X \in \n_\la, X \perp Z,\\
-\tfrac12 \la^2 X - \tfrac14 J_Z^2 X & X \in \n_{\la/2},\\
-\tfrac13 \la^2 X & X \in \n_{2\la/3},\\
-\tfrac13 \la^2 X & X \in \n_{\la/3} \cap \Ker J_Z,\\
0 & X \in \n_{\la/3} \cap (\Ker J_Z)^\perp,
\end{cases}
$$
and so 
$\Tr R_Z^2 = \la^4 + \la^4 (n_\la -1) + \tfrac{\la^4}4 n_{\la/2} - \tfrac{\la^2}4 \Tr (J_Z^{}J_Z^t)_{\n_{\la/2}} +
\tfrac{1}{16} \Tr (J_Z^{}J_Z^t)^2_{\n_{\la/2}}
+ \tfrac{\la^4}{9} n_{2 \la/3} + \tfrac{\la^4}{9} n_{2 \la/3}$, using the fact that $i_{\la/3} = n_{2\la/3}$ from \eqref{eq:ila}.
Equating $\Tr R_Z^2$ to the right hand side of \eqref{eq:TrR2} and substituting $\Tr (J_Z^{}J_Z^t)_{|\n_{\la/2}} = \la^2 n_{\la/2}$
from \eqref{eq:la/2} we find $\Tr (J_Z^{}J_Z^t)^2_{|\n_{\la/2}} = \la^4 n_{\la/2}$. As $\dim \n_{\la/2} = n_{\la/2}$ by definition,
and the operator $J_Z^{}J_Z^t$ is symmetric, the equations
$\Tr (J_Z^{}J_Z^t)_{|\n_{\la/2}} = \la^2 n_{\la/2}, \; \Tr (J_Z^{}J_Z^t)^2_{|\n_{\la/2}} = \la^4 n_{\la/2}$ imply that
$(J_Z^{}J_Z^t)_{|\n_{\la/2}} = \la^2 \, \id_{|\n_{\la/2}}$.

2. We again use the second Ledger formula \eqref{eq:ledger}, $\Tr R_X^2 = \Tr R_A^2$, this time with a unit
vector $X \in \n_{\la/2}$. A direct computation based on \eqref{eq:R} and the fact that
$\<J_VX, J_ZX\> = \la^2\, \<V, Z\> \|X\|^2$ for $V, Z \in \n_\la, \; X \in \n_{\la/2}$, which follows from assertion 1, gives
$$
R_XY=
\begin{cases}
-\tfrac14 \la^2 A & Y = A, \\
-\tfrac14 \la^2 Y & Y \in \n_\la,\\
-\tfrac14 \la^2 Y - \tfrac34 J_{[X, Y]} X & Y \in \n_{\la/2}, \; Y \perp X,\\
-\tfrac13 \la^2 Y & Y \in \n_{2\la/3},\\
-\tfrac16 \la^2 Y & Y \in \n_{\la/3}.
\end{cases}
$$
We first compute the trace of $R_X^2$ on $\n_{\la/2}$. Let $Y_1, \ldots, Y_{n_{\la/2}-1}, Y_{n_{\la/2}}=X$ be an orthonormal
basis for $\n_{\la/2}$. Then
\begin{align*}
\Tr (R_X^2)_{|\n_{\la/2}} &= \frac{\la^4}{16} (n_{\la/2} - 1) +
\frac{3 \la^2}{8} \sum_{s=1}^{n_{\la/2}-1} \|[X, Y_s]\|^2 + \frac{9}{16} \sum_{s=1}^{n_{\la/2}-1} \|J_{[X, Y_s]} X\|^2 \\
&= \frac{\la^4}{16} (n_{\la/2} - 1) +
\frac{15 \la^2}{16} \sum_{s=1}^{n_{\la/2}-1} \|[X, Y_s]\|^2,
\end{align*}
as $\|J_{[X, Y_s]} X\|^2 = \la^2 \|[X, Y_s]\|^2$ by assertion 1.
Now, with $Z_1, \ldots, Z_{n_\la}$ an orthonormal basis for $\n_\la$, we have
$$
\sum_{s=1}^{n_{\la/2}-1} \|[X, Y_s]\|^2 = \sum_{s=1}^{n_{\la/2}-1} \sum_{i=1}^{n_{\la}} \<J_{Z_i}X, Y_s \>^2 =
\sum_{i=1}^{n_{\la}} \|J_{Z_i}X\|^2 = \la^2 n_\la.
$$
So $\Tr (R_X^2)_{|\n_{\la/2}} = \tfrac{\la^4}{16} (n_{\la/2} - 1) + \tfrac{15 \la^4}{16} n_\la$, hence
$$
\Tr R_X^2 = \tfrac{\la^4}{16} n_{\la/2} + \la^4 n_\la + \tfrac{\la^4}{6} n_{2\la/3}.
$$
Comparing this to \eqref{eq:TrR2} we get $n_{2\la/3}=0$, and then $n_{\la/3}=0$ by assertion 1 of Lemma~\ref{l:rootspace}.
\end{proof}

As it follows from Lemma~\ref{l:la/2}, if $\la/2 \in \Delta$, the solvmanifold $G$ is a Damek-Ricci
space. Also, if $\Delta = \{\la\}, \; G$ is the hyperbolic space. To finish the proof it remains to show
that the case $\Delta = \{\tfrac13 \la, \tfrac23 \la, \la\}$ is impossible.

\begin{lemma}\label{l:notla/3}
$\Delta \ne \{\tfrac13 \la, \tfrac23 \la, \la\}$.
\end{lemma}
\begin{proof}
Assume $\Delta = \{\tfrac13 \la, \tfrac23 \la, \la\}$.
Combining the results of Lemma~\ref{l:Jacobi} and Lemma~\ref{l:rootspace}, we obtain that
the volume density function $V(t,\phi)$ along the geodesic $\G_\phi$ is
\begin{equation*}
    V(t, \phi) = \la^{-1} \sinh \la t \; (x_0(t, \phi))^{n_\la-1} \, (x(t, \phi))^{n_{2\la/3}} \, (\det v(t,\phi))^{n_{2\la/3}},
\end{equation*}
with the functions $x_0(t, \phi),\; x(t, \phi)$ and the $2 \times 2$-matrix $v(t,\phi)$ determined by
\begin{align*}
    &\ddot{x}_0 = \la^2 x_0, & &x_0 (0, \phi) = 0, & &\dot x_0 (0, \phi) = 1, \\
    &\ddot x = \tfrac{\la^2}{9} (1 + 2 \Phi^2) x, & &x (0, \phi) = 0, & &\dot x (0, \phi) = 1, \\
    &\ddot v  - \frac{2 \la }{\sqrt 3} \, \Phi \, J \dot v
    - \frac {\la^2}{9} \left(\begin{array}{cc} 1+2 \Phi^2 & 4 \sqrt 3 \, q \Phi \\ -2 \sqrt3 \, q \Phi & 4+2 \Phi^2 \\ \end{array}\right) v = 0, &
    &v (0, \phi) = 0, & &\dot v (0, \phi)  = I_2,
\end{align*}
where $J= (\begin{smallmatrix} 0 & 1 \\ -1 & 0 \end{smallmatrix})$.
Clearly, $x_0 = \la^{-1} \sinh \la t$ and
$$
V(t, 0) = \Bigl(\tfrac{\sinh \la t}{\la}\Bigr)^{n_\la} \Bigl(\tfrac{\sinh (\la t/3)}{\la/3}\Bigr)^{2n_{2\la/3}}
\Bigl(\tfrac{\sinh (2\la t/3)}{2 \la/3}\Bigr)^{n_{2\la/3}},
$$
so to prove the Lemma it remains to show that
$$
x(t, \phi) \det v(t, \phi) \ne \tfrac{27}{2 \la^3} \, \sinh^2 \tfrac{\la}{3}t \; \sinh \tfrac{2 \la}{3} t.
$$
This is indeed the case, as the Taylor expansion of $x(t, \phi) \det v(t, \phi)$ at $t = 0$ is
$$
x(t, \phi) \det v(t, \phi) = t^3 + \tfrac{\la^2}{9} t^5 + \tfrac {2 \la^4}{405} t^7  +
\tfrac{\la^6}{4 \cdot 7 \cdot 3^9} \,(81 - 27 \cos^2(\phi) + 15 \cos^4(\phi) - \cos^6(\phi)) \, t^9
+ o(t^9),
$$
which explicitly depends on $\phi$ (the expansion was obtained using \textsc{Maple}).

\end{proof}

\end{document}